\batchmode
\documentclass [11pt]{article}
\usepackage[latin1]{inputenc} 
\usepackage{amssymb}

\newtheorem{theorem}{Theorem}
\newtheorem{corollary}[theorem]{Corollary}

\newtheorem{proposition}[theorem]{Proposition}
\newtheorem{remark}[theorem]{Remark}
\newtheorem{definition}[theorem]{Definition}
\newtheorem{example}[theorem]{Example}

\title{ Geometric properties of rotation minimizing  vector fields along curves in Riemannian manifolds   }
\author{Fernando ETAYO\footnote{Dept. Mathematics, Statistics and Computation.  University of Cantabria. Avda. de los Castros, s/n, 39071 Santander, SPAIN. e-mail: etayof@unican.es}}
\date{}
\begin{document}
\maketitle

\begin{abstract}
Rotation minimizing (RM) vector fields and frames were introduced by Bishop  as an alternative to the Frenet frame. They are used in CAGD because they can be defined even when the curvature vanishes.  Nevertheless, many other geometric properties have not been studied. In the present paper, RM vector fields along a  curve  immersed into a Riemannian manifold are studied when the ambient manifold is the Euclidean 3-space, the Hyperbolic 3-space and a K\"{a}hler manifold.
\end{abstract}



\bigskip

\section{Introduction}

Rotation minimizing frames (RMF) were introduced by Bishop \cite{B} as an alternative to the  Frenet moving frame along a curve $\gamma$ in $\mathbb{R}^{n}$. The Frenet frame is an orthonormal frame which can be defined for curves in $\mathbb{R}^{n}$, as long as the first $n-1$ derivatives are linear independent. In the classical case $n=3$ the Frenet frame is given by the tangent, the normal and the binormal vectors. Generalizations of Frenet apparatus  to Riemannian manifolds  have been done in the past. In \cite{MR} it is proved that two Frenet curves in the spaces of constant curvature $\mathbb{S}^{n}$ and $\mathbb{H}^{n}$ are congruent if and only if their $n-1$ curvatures are equal, thus generalizing the known result for the Euclidean space $\mathbb{R}^{n}$. Besides, they show that  the converse of this theorem is also true, i.e., Frenet's theorem holds for curves in a connected Riemannian manifold $(M,g)$ if and only if $(M,g)$ is of constant curvature.

An RMF along a curve $\gamma =\gamma (t)$ in $\mathbb{R}^{n}$ is an orthonormal frame defined by the tangent vector and $n-1$ normal vectors $N_{i}$, which do not rotate with respect to the tangent, i.e., $N_{i}'(t)$ is proportional to $\gamma '(t)$. Such a normal vector field along a curve is said to be a rotation minimizing vector field (RM vector field, for short). Any orthonormal basis $\{ \gamma '(t_{0}),N_{1}(t_{0}),\ldots ,  N_{n-1}(t_{0})\}$ at a point $\gamma (t_{0})$ defines a unique RMF along the curve $\gamma$. Thus, such an RFM is  uniquely determined modulo a rotation in $\mathbb{R}^{n-1}$, but it can be defined in any situation of the derivatives of $\gamma$.

 Nowadays, RMF are widely used in Computer Aided Geometric Design (see, e.g., \cite{F}), in order to define a swept surface  by sweeping out a profile in planes normal to the curve. As it is pointed out in \cite{G}, the Frenet frame may result a poor choice for motion planning or swept surface constructions, since it incurs unnecessary rotation of the basis vectors in the normal plane. The fact that the principal normal vector always points to the center of curvature often yields awkward-looking motions, or unreasonably
twisted swept surfaces. Besides, in the points where the curvature vanishes one cannot define the Frenet frame. RM frames avoid these drawbacks, thus being widely used in Computer Aided Geometric Design.  It is a very remarkable fact that Bishop had introduced RM frames before they were interesting in Computer Aided Geometric Design.

\bigskip

In the case of a curve $\gamma$ in an $n$-dimensional Riemannian manifold  $(M,g)$ such an RFM is given (see \cite{A, Et, M}) by a moving orthonormal frame along the curve, $\{ \gamma '(t), N_{1}(t), \ldots ,N_{n-1}(t)\}$, where $\nabla _{\gamma '(t)}N_{i}(t)=-\kappa _{i}(t) \gamma '(t), \quad i=1,\ldots ,n-1$, thus meaning normal vectors $N_{i}$ do not rotate with respect to the tangent vector $\gamma '$. The quantities $\kappa _{i}(t)$ are called the \textit{natural curvatures} and they are functions along the curve. Each of the vectors of the RMF is said to be a rotation minimizing vector. Of course, if $(M,g)$ is the Euclidean space $\mathbb{R}^{n}$, then the notion of RMF particularizes to that of Bishop. This is carefully proved in \cite{Et}.

Let $\nabla$ be the Levi Civita connection of $g$. Then, Frenet type equations read as (see \cite{M, SW})

\begin{equation}
\left( \begin{array}{ccccc} 0 & -\kappa _{1}(t) & -\kappa _{2}(t) & \ldots & -\kappa _{2n-1}(t)\\ \kappa _{1}(t) & 0 & 0 & \ldots & 0\\ \kappa _{2}(t)& 0 & 0&\ldots &0 \\ \vdots & \vdots & \vdots & \ddots & \vdots\\ \kappa _{2n-1}(t) & 0 & 0 & \ldots & 0
\end{array} \right) ,
\label{Frenet}
\end{equation}

\noindent where columms denote the coordinates of the covariant derivatives $\nabla _{\gamma '(t)} \gamma '(t), \quad \nabla _{\gamma '(t)} N_{i}(t)$, $ i=1,\ldots ,n-1 $, of each term of the RMF with respect to this frame.

\bigskip

RM frames in Riemannian manifolds are used in the study of  the structure equations for the evolution of a curve embedded in an n-dimensional Riemannian
manifold with constant curvature  (see, e.g., \cite{M, SW}) or a symmetric Riemannian space (see \cite{A}). They are also used in the study of  mathematical models of equilibrium configurations of thin elastic rods (see, e.g, \cite{Kw}  and the references therein).

\bigskip

The main goal of the present paper is to state geometric properties for RM vector fields along a curve immersed into a Riemannian manifold $(M,g)$. As a formal definition we give the following one:

\begin{definition} Let $\alpha$ be a curve immersed in a Riemannian manifold $(M,g)$.  A normal vector field $N$ over $\alpha$ is said to be an \emph{RM vector field} if it is parallel with respect to
the normal connection of $\alpha$.

\end{definition}

The above condicion is equivalent to the fact $\nabla _{\alpha '}N$ and $\alpha '$ are proportional (see \cite{Et} for the details). As the normal connection is also metric, one can conclude that the norm of an RM vector field is constant and that the angle between two RM vector fields remains constant.

\bigskip

 We focused on three situations, according to the case when the ambient manifold is the Euclidean space, the Hyperbolic space and a K\"{a}hler manifold:

\begin{enumerate}
\item For the case of the Euclidean space $\mathbb{R}^{3}$ we will explicitly show the deep relation between RM vector fields and developable surfaces.

\item In the case of the Hyperbolic space $\mathbb{H}^{3}$ we will show that similar results can be obtained when one has a suitable definition of a developable surface. 

\item For the case of a K\"{a}hler manifold, $J(\gamma ')$ is orthogonal to $\gamma '$, thus  the following question being  natural: is $J(\gamma ')$ always an RM vector field along $\gamma$? Or $\gamma $ is a \textit{special} curve if one can take $N_{1} =J(\gamma ')$, i.e., if $J(\gamma ')$ is an RM vector along $\gamma$? As we will show the answer leads to magnetic curves, which are the integral curves of a convenient 2-form defined by means of the K\"{a}hler form of the manifold.
\end{enumerate}

Finally, we want to point out that some results in the Minkowski space $\mathbb{E}^{n}_{1}$ have been recently obtained for several authors (see, e.g., \cite{Ka}). These are out of the purpose of the present paper.

\section{RM vector fields along curves in $\mathbb{R}^{3}$}

Bishop \cite{B} introduced an RM vector field $N$ over a curve $\alpha$ as a normal vector field along the curve satisfying $N'$ and $\alpha '$ are proportional. In \cite{Et} we have explicitily shown that definition of RM vector field along a curve immersed in a Riemannian manifold  extends that of Bishop:

\begin{theorem} {\rm \cite[Theorem 1]{Et}} A normal vector field $N$ over a curve $\alpha$ immersed in $\mathbb{R}^{3}$ is an RM vector field in the sense of Bishop if and only if it is parallel with 
respect to the normal connection of $\alpha$.
\end{theorem}

The following properties are easy to be proved:

\begin{proposition} Let $\alpha ,\beta$ be two curves immersed in the Euclidean space $\mathbb{R}^{3}$.

\begin{enumerate}
\item The ruled surface defined by a normal  vector field along a curve is developable if and only if the vector field is an RM vector field.

\item If $\alpha$ is the evolute of a curve $\beta$ (and $\beta$ the involute of $\alpha$), then $N(s)=\frac{\beta (s)-\alpha (s)}{\| \beta (s)-\alpha (s) \|}$ defines an RM vector field along $\beta$.

\item The ruled surface defined by an RM  vector field along a curve $\alpha$ is a tangential surface.

\end{enumerate}

\end{proposition}

\textit{Proof.}

\begin{enumerate}
\item The ruled surface can be parametrized as $f(s,\lambda )=\alpha (s)+\lambda N(s)$, with $\alpha$ a unit speed curve and $\| N(s)\| =1$. If $N$ is an RM vector field along $\alpha$, then $[\alpha ',N, N'] =0$, thus proving the surface is developable. If the surface is developable, one has $[\alpha ',N, N'] =0$ and then one can write $N'= a \alpha '+bN$. Taking into account $\| N(s)\| =1$ one obtains $0=N\cdot N'=b$, thus proving $N$ is RM.

\item As is well known, if $\alpha =\alpha (s)$ is a unit speed parametrization of the evolute then $\beta (s)=\alpha (s) +(c-s)\alpha '(s)$ is a parametrization of any involute $\beta$, where $c$ is a constant. A direct calculation shows that $N'(s)$ and $\beta '(s)$ are proportional, thus proving $N$ is an RM vector field along $\beta$.

\item By item 1, that surface  is developable, and then, locally isometric to the plane. Let $f$ be the local isometry.  The locus $\beta$ of the centres of curvature of the curve $f(\alpha )$ is an evolute  of $f(\alpha )$. Then, aplying the inverse local isometry $f^{-1}$ which preserves angles, the given curve $\alpha$ is an involute of $f^{-1}(\beta )$, and the tangential surface to this curve coincides with the given one.
\end{enumerate}

The proof is finished.$\Box$

\bigskip

Item 2 of the above Proposition gives a way to define a Rotation Minimizing Frame (RMF) along a curve $\alpha$, because any curve has infinite evolutes (see, e.g., \cite{E}). Then one can define the RMF given by $\{ \alpha ', N, \alpha ' \times N\}$, where $\times$ denote the cross product in $\mathbb{R}^{3}$.

\bigskip
 The curve in the plane $\mathbb{R}^{2}$ defined by the natural  curvatures $\kappa _{1}, \kappa _{2}$ is said to be the \textit{normal development of the curve} (see \cite{B}). Spherical curves can be characterized by means of their normal development:

\begin{proposition} [\cite{B}] A curve in $\mathbb{R}^{3}$ is spherical if and only if its normal development lies on a line not passing through the origin. The distance of this line from the origin and the radius of the sphere are reciprocals. 

\end{proposition}

The relation between the pair curvature-torsion  $(\kappa ,\tau )$  and the pair of functions $(\kappa _{1} , \kappa _{2})$ is given in the following

\begin{proposition} \emph{\cite[page 52]{McC}} The following relations hold:

$$\kappa =\sqrt{\kappa _{1}^{2}+\kappa _{2}^{2}} \hspace{5mm} {\rm and} \hspace{5mm} \tau = \theta ' =
 \frac{\kappa _{1}\kappa ' _{2}-\kappa ' _{1}\kappa _{2}}{\kappa _{1}^{2} +\kappa _{2}^{2}} \quad ,
$$

\noindent where $\theta =arg(\kappa _{1},\kappa _{2}) =\arctan \frac{\kappa _{2}}{\kappa _{1}}$ and $\theta '$
is the derivative of $\theta$ with respect to the arc length.
\label{natural}

\end{proposition}

Observe that the  normal development of a curve lies on a line  passing through the origin if and only if $\theta '=0$, i.e., if and only if the curve is a plane curve. Ruled surfaces  have been studied in \cite{Tu}, by means of an RMF along the base curve. Assuming $\mathbb{R}^{3}$ is endowed with the Lorentz-Minkowski metric, curves that lie on a surface has been recently characterized by means of RM frames in \cite{Sil}.

\section{RM vector fields along curves in $\mathbb{H}^{3}$}

  As is well known Hyperbolic space can be defined axiomatically as a non-Euclidean geometry. Notions of line and plane can be defined in Hyperbolic 3-space, although relative positions of them are different from that of the Euclidean geometry. By using differential-geometric tools one can study the Hyperbolic space. For instance, lines are geodesics. The first consideration one should have in mind is the existence of different models for $\mathbb{H}^{3}$. All of them are isometric and notions will be introduced without reference to a particular model.

\bigskip

The real Hyperbolic 3-space $\mathbb{H}^{3}$ is the unique up to
isometry 3-dimensional complete, simply connected Riemannian manifold with constant sectional curvature -1. Geodesics of this manifold are called hyperbolic lines. Hyperbolic planes are totally geodesic complete 2-manifolds. For instance, if one consider the Poincar\'{e}'s model  of the upper  hyperspace
$\{ (x,y,z)\in \mathbb{R}^{3}, z>0\}$ with the hyperbolic metric

$$g=\frac{1}{z^{2}}\, (dx^{2}+dy^{2}+dz^{2}) \quad ,
$$

\noindent then hyperbolic lines (resp. planes) are semicircles (resp. hemispheres) orthogonal to the horizontal plane $\{ z=0\}$ and vertical lines (resp. vertical planes). (As this model is conformal,  orthogonality is in both Euclidean and Hyperbolic senses).

\bigskip

The \textit{tangent line} of a curve at a point is the hyperbolic line which is tangent to the curve at the point, i.e., it is the geodesic line through the point with derivative equal to the tangent vector of the curve at the point, as in the Euclidean 3-space where the affine tangent line is the geodesic having the same derivative than the curve. The \textit{tangent plane} of a surface at a point is the hyperbolic plane which is tangent to the surface at the point.

As is well known, the exponential map $\mathrm{exp}_{p}:T_{p}\mathbb{H}^{3}\to \mathbb{H}^{3}$ is a global diffeomorphism. The tangent line $\alpha$ at a point $p=\alpha (s)$ is the image under the exponential map of the line generated by the tangent vector $\alpha '(s)$, and the tangent plane to a surface $S$ at $p$ is $\mathrm{exp}_{p}(T_{p}S)$, where  $T_{p}S$ is the  tangent vector plane to the surface at the point $p$. (In the general case, the exponential map does not send vector subspaces onto totally geodesics submanifolds, but this is the case if the manifold is good enough; see \cite{Ch}).

\bigskip

A \textit{ruled surface} (see \cite{P}) is defined by a smooth family of hyperbolic lines touching a curve, which is called the \textit{directrix} of the surface. Such a surface is said to be  \textit{developable} if the tangent plane of the surface at a point coincides with that at any point of the same  line. As in the Euclidean case, one can parametrize a ruled surface as $f(s,\lambda )= \gamma _{N(s)}(\lambda)$, where $\alpha =\alpha (s)$ is the directrix, parametrized as a unit-speed curve if necessary, and $N(s)$ is the unit vector field along $\alpha$ defining the hyperbolic line $\gamma _{N(s)}$ through the point $\alpha (s)$ by the conditions

$$\gamma _{N(s)}(0)=\alpha (s), \quad \gamma '_{N(s)}(0)=N(s).
$$ 

\bigskip

The following result will be essential in our work.

\begin{proposition}{\rm \cite[Theorem 1]{P}}  A ruled surface $f=f(s,\lambda )=  \gamma _{N(s)}(\lambda)$ is developable if and only if the tangent $\alpha '$ to the directrix, the unit vector $N$ giving the direction of the hyperbolic line of the rulling, and the covariant derivative of the latter along the directrix, $\nabla _{\alpha '}N$, are linearly dependent at any point of the directrix.
\label{Portnoy}
\end{proposition}

This result is independent from the model of the hyperbolic 3-space, because all the models are isometric. The proof given by Portnoy in \cite{P} uses the Poincar\'{e}'s model given by the upper half-space.
Developable surfaces are intensively studied in that paper, where it is proved that a developable surface is isometric to the hyperbolic plane and, reciprocally, a surface having the same intrinsic curvature as that of a hyperbolic plane is necessarily developable. In particular, the \textit{tangential surface} defined by a curve is that defined by the tangent lines to the curve. By the above theorem, it is a developable surface.

\bigskip

We introduce the following

\begin{definition} Let $\alpha, \beta$ be two curves immersed in the Hyperbolic space $\mathbb{H}^{3}$. The curve $\alpha$ is said an evolute of $\beta$ and $\beta $ is said an involute of $\alpha$ if $\beta$ is contained into the tangential surface of $\alpha$ and meets orthogonally the tangent lines of $\alpha$.

\end{definition}

Observe that one can parametrizes the tangential surface to $\alpha$ as $f(s,\lambda )=  \gamma _{\alpha '(s)}(\lambda)$, and an involute $\beta$ as $\beta (s)=\gamma _{\alpha '(s)}(\lambda (s))$. We will not need the explicit determination of the function $\lambda = \lambda (s)$.

\bigskip

We can prove the following results, similar to those of the Euclidean case.

\begin{theorem} The ruled surface defined by a normal  vector field along a curve in $\mathbb{H}^{3}$ is developable if and only if the vector field is an RM vector field.
\label{developable}
\end{theorem}

\textit{Proof.} Let $f(s,\lambda )=  \gamma _{N(s)}(\lambda)$ be a parametrization of the ruled surface with directriz $\alpha =\alpha (s)$. 

If $N$ is an RM vector field, then $\nabla_{\alpha '}N$ and $\alpha '$ are proportional, and the result follows directly from Proposition \ref{Portnoy}.

Let us assume the surface is developable. Then at any point of the curve the following vectors are linearly dependent: $\alpha ', N, \nabla _{\alpha '}N$, which allows us to write  
$\nabla _{\alpha '}N = a \alpha '+ b N
$. Taking into account that $N$ is a unit normal vector field one has:

$$g(\nabla_{\alpha '}N,N)=g(a \alpha '+ b N,N)=b.
$$

From the  identity $g( \nabla_{X}Y,Z)+g( Y, \nabla_{X}Z)  =X(g(Y,Z))$, when one consider $X$ a vector extension of $\alpha'$, and $Y=Z$ unit vector extensions of $ N$, one obtains

$$2g(\nabla_{\alpha '}N,N)=
g(\nabla_{\alpha '}N,N)+
g(N,\nabla_{\alpha '}N)=\alpha '(g(N,N))=0,
$$

\noindent which shows $b=0$. Then one has $\nabla _{\alpha '}N = a \alpha '$, thus proving $N$ is an RM vector field.$\Box$

\begin{corollary} Let $\alpha, \beta$ be two curves immersed in the Hyperbolic space $\mathbb{H}^{3}$.  Assume that $\alpha$ is the evolute of a curve $\beta$ (and $\beta$ the involute of $\alpha$), and let 
$f(s,\lambda )=  \gamma _{\alpha '(s)}(\lambda)$ be a parametrization of the tangential surface to $\alpha$, and $\beta (s)=\gamma _{\alpha '(s)}(\lambda (s))$ a parametrization of $\beta$. Then, the vector field

$$N(s)=\gamma '_{\alpha '(s)}(\lambda (s))
$$
is an RM vector field along $\beta$.

\label{evolute}
\end{corollary}

\textit{Proof}. The ruled surface defined by $N$ with directrix $\beta$ coincides with the tangential surface of the curve $\alpha$, which is developable, by Proposition \ref{Portnoy}. Then, by Theorem \ref{evolute}, the vector field $N$ along $\beta$ is RM.$\Box$

\section{RM vector fields along curves in K\"{a}hler manifolds}

Let us assume that $(M,J,g)$ is a $2n$-dimensional K\"{a}hler manifold. Let $\Omega$ denote the K\"{a}hler form defined by $\Omega (X,Y)=g(JX,Y)$. As is well known, $J$ is an isometry moving any vector to a normal one. If $\gamma$ is a curve immersed in such a manifold, then $J(\gamma ')$ is a normal vector field along the curve and it is natural to ask about the conditions which are satisfied by the curve $\gamma$ in order $J(\gamma ')$ to be an RM vector field. We obtain:

\begin{proposition} Let $\gamma = \gamma (t)$ be a curve in a K\"{a}hler manifold.

\begin{enumerate}

\item Then the vector field $J(\gamma ')$ is RM if and only if $\nabla _{\gamma '(t)} \gamma '(t) =\kappa _{1}(t) J(\gamma '(t))$. In this case, if $\kappa _{1}(t)\equiv 0$, then $\gamma $ is a geodesic.
\item If $J(\gamma ')$ is RM then $\nabla _{\gamma '(t)}N_{i}(t)=0, \quad i=2,\ldots ,2n-1$, for all normal vector fields $N_{i}, \quad i=2,\ldots ,2n-1$ such that $\{ \gamma ', J(\gamma '),N_{2} \ldots ,N_{2n-1}\}$ is an RMF, i.e., the natural curvatures $\kappa _{2},\ldots ,\kappa _{2n-1}$ vanish.
\item If $J(\gamma ')$ is RM then  $\parallel \gamma '(t) \parallel= \sqrt{ g(\gamma '(t),\gamma '(t))}$ is constant.

\end{enumerate}
\label{basico}
\end{proposition}

\textit{Proof.}
\begin{enumerate}

\item As $(M,J,g)$ is  K\"{a}hler one has $\nabla J=0$. A direct computation shows $ \nabla _{\gamma '} J(\gamma ') =-\kappa _{1}(t)\gamma '$ if and only if $\nabla _{\gamma '} \gamma ' =\kappa _{1}(t)J(\gamma ')$. If $\kappa _{1}(t)\equiv 0$ then $\nabla _{\gamma '} \gamma ' \equiv 0$, thus proving $\gamma$ is a geodesic.

\item It is a direct  consequence of expression (\ref{Frenet}).

\item Taking into account the properties of the Levi Civita connection $\nabla $ of $g$ one has

\begin{eqnarray*}
0=(\nabla _{\gamma '}g)(\gamma ',\gamma ')=\gamma '(g(\gamma ',\gamma '))-2g(\nabla _{\gamma '}\gamma ',\gamma ') =
\gamma '(g(\gamma ',\gamma '))- 2g(\kappa _{1}J (\gamma '),\gamma ')=\\
\gamma '(g(\gamma ',\gamma '))- 2\kappa _{1}\Omega (\gamma ',\gamma ')  =
\gamma '(g(\gamma ',\gamma ')),
\end{eqnarray*}
thus proving $g(\gamma ',\gamma ')$ is  constant along $\gamma$. $\square$
\end{enumerate}
\bigskip

Remember the following

\begin{definition} {\rm (See \cite{otta, miho} and \cite[p.418] {metal} )}. An analytically  planar curve in a K\"{a}hler manifold $(M,J,g)$ is a curve such that $\nabla _{\gamma '} \gamma ' = a(t)\gamma ' + b(t) J(\gamma ')$, where $a,b$ are functions on the curve.
\end{definition}

The above analytically planar curves are also often called $h$-planar, holomorphically planar, $H$-planar or $J$-planar curves. These curves are special cases of quasigeodesic \cite{pe} and $F$-planar curves \cite{misi} and \cite[p.385]{metal}. A curve having $J(\gamma ')$ as an RM vector field is an analytically  planar curve. Besides, when  $\kappa _{1}$ is constant, the curve  is also a magnetic curve, because of the following 

\begin{definition} {\em (See \cite{A2011})}.
A curve satisfying $\nabla _{\gamma '} \gamma ' =\kappa _{1}J(\gamma ')$ with $\kappa _{1}\in \mathbb{R}$ a real constant, is said to be a \textit{magnetic curve} or a \textit{trajectory of the magnetic field} given by the 2-form $\kappa _{1} \Omega$, where $\Omega$ is  the K\"{a}hler form of $(M,J,g)$.
\label{defAdachi}
\end{definition}

If $J(\gamma ')$ is an RM vector, with $\kappa _{1}(t)=\kappa _{1}$ a real constant, then $\gamma$ is a magnetic curve with respect to the  2-form $\kappa _{1} \Omega$, thus allowing one to apply all the known results for this kind of curves. One has:

\begin{theorem} Let $\gamma$ be a curve in a K\"{a}hler manifold $(M,J,g)$ and let us assume that $J(\gamma ')$ is an RM vector along $\gamma$, such that $ \nabla _{\gamma '} J(\gamma ') =-\kappa _{1}\gamma '$, with $\kappa _{1}$ a real constant. Then one has:

\begin{enumerate}

\item The curve $\gamma$  is a  magnetic curve with respect to  the 2-form $\kappa _{1} \Omega$.

\item {\rm \cite[Theorem 4]{K}} If $(M,J,g)$ has constant holomorphic curvature, then the curve $\gamma$ is contained in a totally geodesic surface in $M$.
\end{enumerate}

\end{theorem}

Last item of the above theorem agrees with the vanishing of the last  natural curvatures $\kappa _{2},\ldots ,\kappa _{2n-1}$ obtained in Proposition \ref{basico}. At the points of the curve, vectors $\gamma '$ and $J(\gamma ')$ define a basis of the tangent plane  of the totally geodesic surface in which the curve is immersed, and then, as this surface is totally geodesic and $N_{i}=0, \quad i=2,\ldots ,2n-1$, are normal to the surface, one has $\nabla _{\gamma '}N_{i}=0$.

\begin{example} {\rm (See \cite{A1995} and \cite[Examples 1,2 3]{A2011})}. Let $\gamma$ be a curve in a  complex space form such that $J(\gamma ')$ is an RM vector along $\gamma$ with $\kappa _{1}\in \mathbb{R}$ a real constant, and let us assume $\kappa _{1}\neq 0$. Then one has:

\begin{enumerate}
\item If $M=\mathbb{C}^{n}$, then $\gamma$ is  a  circle.

\item If $M=\mathbb{C}P^{n}(c)$, then $\gamma$ is a small circle in some totally  geodesic $\mathbb{C}P^{1} \cong S^{2}$.

\item If $M=\mathbb{C}H^{n}(-c)$, then $\gamma$ is a line in a totally geodesic $\mathbb{C}H^{1}\cong H^{2}$.
\end{enumerate}
\label{exkconstante}
\end{example}

  In a more general context one has the following

\begin{definition} {\em (See \cite{Ba})} A curve $\gamma$ is said to be a trajectory of the magnetic field given by a 2-form $F$ if  $\nabla _{\gamma '} \gamma ' =\Phi (\gamma ')$, where $\Phi$  is the operator defined by the relation $g(\Phi (X),Y)=F(X,Y)$.
\label{defBarros}
\end{definition}

Definition \ref{defAdachi} is a particular case of Definition \ref{defBarros}, taking $\Phi = \kappa _{1} J$ and $F=\kappa _{1}\Omega$.  Obviously,  $J(\gamma ')$ is an RM vector if and only if $\gamma$ is a magnetic curve for $F=f\Omega$, $f$ being any smooth extension of $\kappa _{1}$ to the manifold $M$.

\bigskip

First, we are interested in the case where $\kappa _{1}(t)$ is a non-constant function. Let $\gamma =\gamma (s)$ be a unit speed curve in $\mathbb{C}=\mathbb{R}^{2}$. In this case,  by Formula \ref{Frenet},  the natural curvature $\kappa _{1}=\pm \kappa$ (see also Proposition \ref{natural}, taking into account that the torsion $\tau =0$).  Then $J(\gamma ')$ is an RM vector field along $\gamma$ if and only if $(J(\gamma '))'= -\kappa _{1} \gamma '$. A direct calculation shows that $J(\gamma ')$ is an RM vector field along $\gamma$ if and only the following system of differential equations

\begin{equation}
\left\{ \begin{array}{c}\gamma _{2}''(s) =\kappa _{1} (s)\gamma _{1}'(s)\\
\gamma _{1}''(s)= -\kappa _{1}(s) \gamma _{2}'(s)\end{array}
\right\}
\label{ecuaciones}
\end{equation}
is satisfied, defining the complex structure $J$ as usual by

$$
J \left( \begin{array}{c} a\\ b \end{array} \right)=
\left( \begin{array}{cc}
0 & -1 \\ 1 & 0 \end{array} \right)
\left( \begin{array}{c} a\\ b \end{array} \right)=
\left( \begin{array}{c} -b\\ a \end{array} \right) .
$$ 

System \ref{ecuaciones} can be found in any book of Differential Geometry when Frenet equations are integrated in the case of a plane curve (see, for instance \cite{E}). Thus, one cannot go forward: the problem of finding curves in $\mathbb{C}$ having $J(\gamma ')$ as an RM vector field is equivalent to that of finding a unit speed parametrization of the curve.

If $\kappa _{1}(t)=\kappa _{1}$ is constant one can solve explicitly the system, obtaining:

\begin{equation}
\left\{ \begin{array}{c}
\gamma _{1}(t) =A_{1} + B\cos (-\kappa _{1}t) +C\sin (-\kappa _{1}t) \\
\gamma _{2}(t)= A_{2} - B\sin (-\kappa _{1}t) +C\cos (-\kappa _{1}t) \end{array} 
\right\} ,
\label{soluciones}
\end{equation}
which are circles with center $(A_{1},A_{2})$ and radius $\sqrt{B^{2}+C^{2}}$.

Dividing both equations in (\ref{ecuaciones}), one also can solve the system in the general case of $\kappa _{1}(t)$ being a function with $\kappa _{1}(t)\neq 0, \forall t$. One obtains

\begin{equation}
0=\gamma _{1}'(t)  \gamma _{1}''(t) + \gamma _{2}'(t) \gamma _{2}''(t) =\frac{1}{2} ( (\gamma _{1}'(t))^{2}+(\gamma _{2}'(t))^{2})'
=\frac{1}{2} (\parallel \gamma '(t)\parallel^{2})' \quad ,
\label{norma}
\end{equation}
thus proving the norm is constant. Besides, in this case, Equations (\ref{ecuaciones}) and (\ref{norma}) are equivalent, thus proving any curve of constant speed  has $J(\gamma ')$ as an RM vector field  (by Proposition \ref{basico}, item 3 we knew one of the implications: if $J(\gamma ')$ is an RM vector field then $\parallel \gamma '\parallel$ is constant). As any curve has a natural parametrization, one can always re-parametrize the curve to satisfy equation (\ref{norma}). 

\begin{example}
\label{espiral}
For instance, consider the logarithmic spiral $\gamma (t)=(e^{t}\cos t,e^{t}\sin t)$. A natural parametrization for this curve is given by

$$ \gamma (s)=\left( (1+\frac{s}{\sqrt{2}})\cos (\log (1+\frac{s}{\sqrt{2}})) \quad,\quad (1+\frac{s}{\sqrt{2}}) \sin (\log (1+\frac{s}{\sqrt{2}}))\right), \quad s>0.
$$
A direct computation shows that $(J(\gamma '))'= -\kappa _{1} \gamma '$ with $\kappa _{1}(s)=  (-1)/(s+\sqrt{2})$. It is easily shown that $(J(\gamma '))'= -\kappa _{1} \gamma '$ has no solution for $\gamma =\gamma (t)$.
\end{example}

\begin{remark}
The situation  can be generalized to any Riemannian surface $(M,J,g)$ in the sense that any curve of constant speed has  $J(\gamma ')$ as an RM vector field (see \cite{Ba}). 
\end{remark}

\bigskip

\begin{remark} Let $\gamma$ be a curve in $\mathbb{C}^{2}$. Then $J(\gamma ')$ is an RM vector field if and only if the following system of ODE 

\begin{equation}
\left\{ \begin{array}{c}\gamma _{3}''(t) =\kappa _{1}(t) \gamma _{1}'(t)\\
\gamma _{4}''(t)= \kappa _{1}(t) \gamma _{2}'(t)\\
\gamma _{1}''(t) =-\kappa _{1} (t)\gamma _{3}'(t)\\
\gamma _{2}''(t)= -\kappa _{1}(t) \gamma _{4}'(t)
\end{array}
\right.
\label{ecuacionesc2}
\end{equation}
is satisfied. When we are working in complex dimensions greater than one, not any constant-speed curve has $J(\gamma ')$ as an RM vector field. For example, consider the curve 

$$\gamma (s) =(\cos s,\sin s ,0,0)$$ 

\noindent in $\mathbb{C}^{2}$. In this case, $J(\gamma ')$ is not an RM vector field. In fact, any solution of Equation (\ref{ecuacionesc2}) with $\kappa _{1}$ a non zero constant, is a circle, as we have said in Example \ref{exkconstante}, but not any circle has the property of $J(\gamma ')$ being an RM vector field.

\end{remark}

\bigskip

\textbf{Acknowledgements.} The author wants to thank his colleagues L. A. Fern\'{a}ndez (Universidad de Cantabria) and Luiz C. B. da Silva (Universidade Federal de Pernambuco) for their useful comments. Besides, the author acknowldges the suggestions and comments written by the referees, which have allowed to improve the paper.

\end{document}